\def\version{08/08/2008\quad Version 5}
\documentclass[11pt,a4paper]{amsart}
\usepackage{fullpage}
\usepackage{amssymb,amscd}
\usepackage[all]{xy}
\usepackage{mathrsfs}
\usepackage{hyperref}
\usepackage[numeric,lite]{amsrefs}

\newtheorem{thm}{Theorem}[section]
\newtheorem{theorem}[thm]{Theorem}
\newtheorem{lemma}[thm]{Lemma}
\newtheorem{proposition}[thm]{Proposition}
\newtheorem{corollary}[thm]{Corollary}

\theoremstyle{definition}
\newtheorem{remark}[thm]{Remark}
\newtheorem{definition}[thm]{Definition}

\numberwithin{equation}{section}
\numberwithin{figure}{section}

\def\ie{\emph{i.e.}}

\def\:{\colon}
\def\.{\cdot}
\def\o{\circ}
\def\<{\left\langle}
\def\>{\right\rangle}
\def\({\left(}
\def\){\right)}
\def\ph#1{\phantom{#1}}
\def\epsilon{\varepsilon}
\def\phi{\varphi}

\def\leq{\leqslant}
\def\geq{\geqslant}

\def\lra{\longrightarrow}

\def\ra{\rightarrow}
\def\bar#1{\overline{#1}}

\def\iso{\cong}

\DeclareMathOperator{\aug}{\epsilon}
\DeclareMathOperator{\id}{id}

\DeclareMathOperator{\im}{im}
\DeclareMathOperator{\rank}{rank}

\def\CP{\mathbb{C}\mathrm{P}}

\def\F{\mathbb{F}}

\def\k{\Bbbk}

\def\Z{\mathbb{Z}}

\DeclareMathOperator{\cone}{C}

\DeclareMathOperator{\HAQ}{HAQ}

\DeclareMathOperator{\TAQ}{TAQ}
\DeclareMathOperator{\Tor}{Tor}

\newdir{ >}{{}*!/-8pt/@{>}}

\title[Topological Andr\'e-Quillen homology for cellular $S$-algebras]
{Topological Andr\'e-Quillen homology for cellular commutative
\boldmath$S$-algebras}
\date{\version \hfill\emph{arXiv:0708.2041}}
\author{Andrew Baker \and Helen Gilmour \and Philipp Reinhard}
\address{
Department of Mathematics, University of Oslo, Norway.
\hfill\newline
Permanent address: Department of Mathematics, University of Glasgow,
Glasgow G12 8QW, Scotland.}
\email{a.baker@maths.gla.ac.uk}
\urladdr{http://www.maths.gla.ac.uk/$\sim$ajb}
\email{h.s.e.gilmour@durham.ac.uk}
\email{p.reinhard@maths.gla.ac.uk}
\thanks{
A.~Baker was partially supported by a YFF Norwegian Research Council
grant while at the University of Oslo in 2007--8, and Intas grants
03-51-3251 and 06-1000017-8609; H.~Gilmour was supported by an EPSRC
studentship; P.~Reinhard was supported by an ORS grant. We would like
to thank M.~Basterra, P.~Kropholler, M.~Mandell, P.~May, B.~Richter,
J.~Rognes and S.~Sagave for numerous helpful comments. We are also very
grateful to the referee for encouraging us to rethink significantly
issues of notation and structure, thus improving the structure of the
paper.
}

\begin{document}

\begin{abstract}
Topological Andr\'e-Quillen homology for commutative $S$-algebras
was introduced by Basterra following work of Kriz, and has been
intensively studied by several authors. In this paper we discuss
it as a homology theory on CW commutative $S$-algebras and apply
it to obtain results on minimal atomic $p$-local $S$-algebras which
generalise those of Baker and May for $p$-local spectra and simply
connected spaces. We exhibit some new examples of minimal atomic
commutative $S$-algebras.
\end{abstract}

\maketitle

\section*{Introduction}

In this paper we give an account of some results on \emph{topological
Andr\'e-Quillen homology} and \emph{cohomology} for CW commutative
$A$-algebras, where $A$ is a commutative $S$-algebra. The main goal
is to develop arguments based on skeletal filtrations with a view
to continuing the work begun in~\cite{Hu-Kriz-May} by extending
results of~\cite{AJB-JPM} to the case of CW~commutative $S$-algebras.
Some of this work originally appeared in the second author's PhD
thesis~\cite{HG:PhD}, but we go further and use it to investigate
some examples. Our main sources on topological Andr\'e-Quillen
(co)homology include
\cites{IK:BP,Basterra-TAQ,Basterra-Mandell,Basterra-McCarthy,Basterra-Richter,Mandell-TAQ}.

For definiteness, we work in the model category of commutative
$S$-algebras described in~\cite{EKMM}; following the referee's
suggestions we write cofibration or fibration in place of
$q$-cofibration or $q$-fibration. We use the standard notions
$\twoheadrightarrow$, $\rightarrowtail$ and $\xrightarrow{\sim}$
to denote fibrations, cofibrations and weak equivalences, respectively.
For general notions of model categories, see~\cite{Hovey:ModCats}.
Given a model category $\mathscr{M}$, we will denote its homotopy
category by $\bar{h}\mathscr{M}$. Thus for a commutative $S$-algebra
$A$, $\mathscr{M}_A$ denotes the category of $A$-modules and
$\bar{h}\mathscr{M}_A$ denotes its homotopy category; the latter
is written $\mathscr{D}_A$ in~\cite{EKMM}, but we
follow~\cite{Basterra-TAQ}. Where necessary, we will assume that
(co)fibrant replacements are made.

We will often consider a map of commutative $S$-algebras $A\lra B$.
The notations $B/A$ and $B|A$ are used in Basterra~\cite{Basterra-TAQ}
and Quillen~\cite{DQ:CoHomRings} to indicate such a pair of $S$-algebras.
This notation is not ideal given the appearance of algebras over and
under a given one, therefore we follow the suggestions of the referee
in adopting alternatives which we hope are more suitable. In particular,
we use the traditional `pair' notation when discussing (co)homology,
writing $(B,A)$ for $A\lra B$, thus following~\cite{Basterra-Mandell}.

\section{Recollections on Topological Andr\'e-Quillen (co)homology}
\label{sec:TAQhomology-recall}

Throughout, we always assume that $A$ be a cofibrant commutative
$S$-algebra. We recall from~\cite{Basterra-TAQ} that for a pair of
commutative $S$-algebras $(B,A)$, there is a $B$-module $\Omega_A(B)$
(well defind in the homotopy category $\bar{h}\mathscr{M}_B$) for
which
\begin{equation}\label{eqn:Omega-char}
\bar{h}\mathscr{C}_{A}/B(B,B\vee M)\iso\bar{h}\mathscr{M}_B(\Omega_A(B),M).
\end{equation}
Here $\bar{h}\mathscr{C}_{A}/B$ denotes the derived category of commutative
$A$-algebras over $B$. Notice that the identity map on $\Omega_A(B)$
corresponds to a morphism $B\lra B\vee\Omega_A(B)$ which projects onto
the \emph{universal derivation}
\begin{equation}\label{eqn:Omega-UnivDer}
\delta_{(B,A)}\in\bar{h}\mathscr{M}_A(B,\Omega_A(B)).
\end{equation}

The \emph{topological Andr\'e-Quillen homology} and \emph{cohomology}
of $(B,A)$ with coefficients in a $B$-module $M$ are defined by
\begin{subequations}\label{eqn:TAQ*}
\begin{align}
\TAQ_*(B,A;M)&=\pi_*\Omega_A(B)\wedge_BM,
\label{eqn:TAQ_*}\\
\TAQ^*(B,A;M)&=\pi_{-*}F_B(\Omega_A(B),M),
\label{eqn:TAQ^*}
\end{align}
\end{subequations}
where $F_B$ denotes the internal function object in $\mathscr{M}_B$.

Associated to an $A$-algebra map $B\lra C$, there are natural long
exact sequences
\begin{subequations}\label{eqn:TAQ-LES}
\begin{multline}\label{eqn:TAQ-LES_*}
\cdots\lra\TAQ_k(B,A;M)\lra\TAQ_k(C,A;M)\lra\TAQ_k(C,B;M) \\
\lra\TAQ_{k-1}(B,A;M)\lra\cdots
\end{multline}
\begin{multline}\label{eqn:TAQ-LES^*}
\cdots\lra\TAQ^k(C,B;M)\lra\TAQ^k(C,A;M)\lra\TAQ^k(B,A;M) \\
\lra\TAQ^{k+1}(C,B;M)\lra\cdots
\end{multline}
\end{subequations}

Given a pair $(B,A)$, suppose that $E$ is a $B$ ring spectrum with
product $\mu\:E\wedge_B E\lra E$, unit $j\:B\lra E$ and induced unit
over $A$, $j_0\:A\lra B\xrightarrow{j} E$. Then suppressing mention
of natural isomorphisms of $A$-modules of form $A\wedge_A X\iso X$,
we obtain a commutative diagram in the homotopy category
$\bar{h}\mathscr{M}_B$.
\begin{equation}\label{eqn:theta-defndiagram}
\xymatrix{
& \Omega_{A}(B)\ar[r]^{j_0\wedge\mathrm{id}}\ar[ddd]
& E\wedge_A\Omega_{A}(B)\ar[dr]\ar[ddd] &  \\
B \ar[ur]^{\delta_{(B,A)}}\ar[d]_{j_0}\ar@{.}'[r][rr]
& &\ar@{.>}[r] & E\wedge_{B}\Omega_{A}(B) \\
E\wedge_A B\ar[dr]_{\mathrm{id}\wedge\delta_{(B,A)}}
& & & E\wedge_B E\wedge_B\Omega_{A}(B)
\ar[u]_{\mu\wedge\mathrm{id}} \\
& E\wedge_A \Omega_A(B) \ar[r]_{\mathrm{id}\wedge j_0\wedge\mathrm{id}\ph{-}}
& E\wedge_A E\wedge_A\Omega_A(B)\ar[ur] &
}
\end{equation}
Now applying homotopy to the dotted line, \ie, the composition
\begin{equation}\label{eqn:Hurewicz-UnivDeriv-E}
B \xrightarrow{\delta_{(B,A)}} \Omega_A(B) = A\wedge_A\Omega_A(B)
\lra E\wedge_A\Omega_A(B)\lra E\wedge_B\Omega_A(B),
\end{equation}
we obtain the \emph{$\TAQ$ Hurewicz homomorphism}
\begin{equation}\label{eqn:theta-defn}
\theta\:\pi_*B \lra \TAQ_*(B,A;E).
\end{equation}

Now let $\phi\:B\lra C$ be a morphism of commutative $A$-algebras
for which $E$ also admits the structure of a $C$ ring spectrum so
that $\phi^*\:\bar{h}\mathscr{M}_C\lra\bar{h}\mathscr{M}_B$ gives
the same $B$ ring spectrum structure. Then the following diagram
commutes, where the morphisms are the evident ones.
\begin{equation}\label{eqn:Hurewicz-diagram}
\xymatrix{
B\ar[r]^{\ph{A}\delta_{(B,A)}\ph{AB}}\ar[d] & \Omega_A(B)\ar[r]\ar[d]
                                  & E\wedge_B\Omega_A(B)\ar[d] \\
C\ar[r]^{\ph{A}\delta_{(C,A)}\ph{AB}}\ar[dr]_{\delta_{(C,B)}}
& \Omega_A(C)\ar[r]\ar[d] & E\wedge_C\Omega_A(C)\ar[d] \\
& \Omega_B(C)\ar[r] & E\wedge_C\Omega_B(C)
}
\end{equation}
This shows that the Hurewicz homomorphism is natural with respect
to morphisms of pairs $\phi\:(B,A)\lra(C,A)$ and $C$ ring spectra
$E$ in the sense that the diagram
\begin{equation}\label{eqn:Hurewicz-natural}
\xymatrix{
\pi_*B \ar[r]^{\theta\ph{----}} \ar[d]_{\phi_*}
               & \TAQ_*(B,A;E) \ar[d]^{\phi_*} \\
\pi_*C \ar[r]^{\theta\ph{----}} & \TAQ_*(C,A;E)
}
\end{equation}
commutes. In our work we will use the Hurewicz homomorphism
$\theta\:\pi_*B \lra \TAQ_*(B,A;E)$ when $E$ is also a commutative
$B$-algebra, and use the naturality when there are morphisms of
$A$-algebras $B\lra C\lra E$.

We will also use the following result which is immediate from the
definition. Here we use the notation $E_*^A(-)=\pi_*E\wedge_A(-)$
and the Hurewicz homomorphism $\underline{e}\:\pi_*B\lra E_*^A B$
is induced by the unit of $E$ over $A$.
\begin{proposition}\label{prop:Hurewicz-E}
The $\TAQ$ Hurewicz homomorphism factors as
\[
\theta\:\pi_*B \xrightarrow{\;\underline{e}\;} E_*^A(B)
       \xrightarrow{\;\ph{\underline{e}}\;}\TAQ_*(B,A;E).
\]
\end{proposition}

We will be especially interested in the situation where
$A$ and $B$ are connective and the map $\phi\:A\lra B$
induces an isomorphism $\pi_0A\xrightarrow{\iso}\pi_0B$;
we will write $\k=\pi_0A=\pi_0B$. Then there is an
Eilenberg-Mac~Lane object $H\k$, which can be taken to
be a CW commutative $A$-algebra or $B$-algebra, which
allows us to define the \emph{ordinary topological
Andr\'e-Quillen homology} and \emph{cohomology} of $(B,A)$:
\begin{subequations}\label{eqn:HAQ*}
\begin{align}
\HAQ_*(B,A)&=\TAQ_*(B,A;H\k)=\pi_*\Omega_A(B)\wedge_BH\k,
\label{eqn:HAQ_*}\\
\HAQ^*(B,A)&=\TAQ^*(B,A;H\k)=\pi_{-*}F_B(\Omega_A(B),H\k).
\label{eqn:HAQ^*}
\end{align}
\end{subequations}
When $\pi_0C=\k$, the long exact sequences of~\eqref{eqn:TAQ-LES}
yield long exact sequences in $\HAQ_*$ and $\HAQ^*$. We can
also introduce coefficients in a $\k$-module $M$ by setting
\begin{subequations}\label{eqn:HAQ*-Coeffs}
\begin{align}
\HAQ_*(B,A;M)&=\TAQ_*(B,A;HM)=\pi_*\Omega_A(B)\wedge_BHM,
\label{eqn:HAQ_*-Coeffs}\\
\HAQ^*(B,A;M)&=\TAQ^*(B,A;HM)=\pi_{-*}F_B(\Omega_A(B),HM).
\label{eqn:HAQ^*-Coeffs}
\end{align}
\end{subequations}

An important result on $\HAQ_*$ is provided by~\cite{Basterra-TAQ}*{lemma~8.2}.
However it appears that this result is incorrectly stated (although
the proof seems to be correct) and should read as follows. For a map
of $A$-modules $\theta$, we denote by $\cone_\theta$ the mapping cone
of $\theta$ in $\mathscr{M}_A$.
\begin{lemma}[Basterra~\cite{Basterra-TAQ}*{lemma~8.2}]\label{lem:Basterra8.2}
Let $\phi\:A\lra B$ be an $n$-equivalence, where $A$ and $B$ are connective
and $n\geq1$. Then $\Omega_A(B)$ is $n$-connected and there is a map of
$A$-modules $\tau\:\cone_\phi\lra\Omega_A(B)$ for which
\[
\tau_*\:\pi_{n+1}\cone_\phi\xrightarrow{\iso}\pi_{n+1}\Omega_A(B).
\]
\end{lemma}

An immediate consequence is an analogue of the classical Hurewicz
theorem.
\begin{corollary}\label{cor:Basterra8.2}
The map $\tau$ induces isomorphisms
\[
\tau_*\:\pi_k\cone_\phi\xrightarrow{\iso}\HAQ_k(B,A)\quad(k\leq n+1).
\]
\end{corollary}
\begin{proof}
{}From~\cite{EKMM} there is a K\"unneth spectral sequence for which
\[
\mathrm{E}^2_{p,q}=\Tor^{B_*}_{p,q}(\pi_*(\Omega_A(B)),\k)
\Longrightarrow
\pi_{p+q}\Omega_A(B)\wedge_{B}H\k=\HAQ_{p+q}(B,A).
\]
For dimensional reasons we have
$\mathrm{E}^\infty_{0,n+1}=\mathrm{E}^2_{0,n+1}$ and so, on recalling
that $\k=\pi_0A=\pi_0B$,
\[
\HAQ_{n+1}(B,A)=[\pi_*(\Omega_A(B))\otimes_{B_*}\k]_{n+1}
                     =\pi_{n+1}\Omega_A(B)\otimes_{B_0}\k
                     =\pi_{n+1}\Omega_A(B).
\qedhere
\]
\end{proof}

Recall that for any $A$-module $X$, there is a free commutative
$A$-algebra on~$X$,
\[
\mathbb{P}_AX=\bigvee_{i\geq0}X^{(i)}/\Sigma_i,
\]
where $X^{(i)}$ is the $i$-th smash power over $A$. We remark
that if $A\lra A'$ is a morphism of commutative $S$-algebras,
then from~\cite{EKMM} we have
\begin{equation}\label{eqn:P-naturality}
\mathbb{P}_{A'}(A'\wedge_{A}X)\iso A'\wedge_{A}\mathbb{P}_AX.
\end{equation}
The $A$-algebra map $\mathbb{P}_AX\lra\mathbb{P}_A*=A$ induced
by collapsing~$X$ to a point makes $A$ into an
$\mathbb{P}_AX$-algebra, and allows us to view $\mathbb{P}_AX$
as an augmented $A$-algebra. There is a cofibration sequence of
$\mathbb{P}_AX$-modules
\begin{equation}\label{eqn:AX->PX->A}
\mathbb{A}_AX\lra\mathbb{P}_AX\xrightarrow{\ \aug\ }\mathbb{P}_A*=A,
\end{equation}
where
\[
\mathbb{A}_AX=\bigvee_{i\geq1}X^{(i)}/\Sigma_i
\]
and $\aug$ is the augmentation. We note that $\mathbb{A}_AX$
is an $A$-nuca, \ie, a non-unital commutative $A$-algebra.
For the $A$-sphere $S^n=S^n_A$ with $n>0$, we obtain the
commutative $A$-algebra $\mathbb{P}_AS^n$ and augmentation
$\mathbb{P}_AS^n\lra A$; this allows us to view an $A$-module
or algebra as a $\mathbb{P}_AS^n$-module or algebra.

For a pair $(B,A)$, the $B$-module $\Omega_A(B)$ is defined
in the homotopy category $\bar{h}\mathscr{M}_B$ by
\begin{equation}\label{eqn:Omega-defn}
\Omega_A(B)=\mathrm{L}Q_B \mathrm{R}I_B(B^c\wedge_A B).
\end{equation}
Here we are using the following notation in keeping with~\cite{Basterra-TAQ}.
Thus $(-)^c$ is a cofibrant replacement functor, while $\mathrm{R}I_B$
is the right derived functor of the augmentation ideal functor on the
category of $B$-algebras augmented over $B$. The targets of $I_B$ and
$\mathrm{R}I_B$ are the category of $B$-nucas (non-unital $B$-algebras)
and its homotopy category. Also, $\mathrm{L}Q_B$ is the left derived
functor of $Q_B$, defined by the following strict pushout diagram in
the category of $B$-modules.
\begin{equation}\label{eqn:Q}
\xymatrix{ N \wedge_B N \ar[d]\ar[r] & {*} \ar[d]\\
N \ar[r] & Q_B(N)}
\end{equation}

\begin{lemma}\label{lem:OmegaPX}
Let $A$ be a commutative $S$-algebra, let $N$ be an $A$-nuca and
let~$B$ be an $A$-algebra. Then
\[
Q_B(B \wedge_A N)=B \wedge_A Q_A(N).
\]
\end{lemma}
\begin{proof}
This follows from the identities
\begin{align*}
B \wedge_A *&  = *, \\
B \wedge_A N \wedge_B B \wedge_A N &\cong B\wedge_A N \wedge_A N,
\end{align*}
together with the fact that the functor $B \wedge_A (-)$ is
a left adjoint in the category of modules, hence it respects
colimits.
\end{proof}

For an $A$-module $M$ and $A$-nuca $N$, the $A$-nuca $\mathbb{A}_AX$
satisfies the adjunction isomorphism
\[
\mathscr{M}_A(X,N) \iso \mathscr{N}_A(\mathbb{A}_AX,N),
\]
where $\mathscr{N}_A$ denotes the category of $A$-nucas. The
idea for the proof of our next result was due to Maria Basterra.
\begin{lemma}\label{lem:Q(AX)}
For each $A$-module $M$, there is a natural isomorphism
\[
\mathscr{M}_A(Q_A(\mathbb{A}_AX),M) \iso \mathscr{M}_A(X,M),
\]
hence $Q_A(\mathbb{A}_AX)\iso X$.
\end{lemma}
\begin{proof}
Given the natural isomorphism, Yoneda's lemma implies that
$Q_A(\mathbb{A}_AX)\iso X$.

Recall Basterra's functor $Z\:\mathscr{M}_A\lra\mathscr{N}_A$
which assigns to each $A$-module $M$ the same module $Z(M)$
with the trivial product $Z(M)\wedge_AZ(M)\lra Z(M)$. Then
for any $A$-nuca $N$ we have
\[
\mathscr{M}_A(Q_A(N),M) \iso \mathscr{N}_A(N,Z(M)),
\]
since an $A$-module map $Q_A(N)\lra M$ fits into a commutative
diagram
\[
\xymatrix{ N \wedge_A N \ar[d]\ar[r] & {*}\ar[d]\ar@/^/[ddr] &\\
N \ar[r]\ar@/_/[drr] & Q_A(N)\ar[dr] & \\
&& M
}
\]
in which the rectangle is a pushout diagram. Hence there is
a factorisation
\[
\xymatrix{
N \wedge_A N \ar[d]\ar[r] & N\ar[d] \\
M \wedge_A M \ar[r]^{\ph{abc}0} & M
}
\]
of the products, while for an $A$-nuca map $N\lra Z(M)$ there
is a factorisation
\[
\xymatrix{
N \wedge_A N \ar[d]\ar[rr] && N\ar[d] \\
Z(M) \wedge_A Z(M) \ar[r] & {*}\ar[r] & Z(M)
}
\]
showing that the map factors through a module map $Q_A(N)\lra M$.

Now for any $A$-module $X$, we have
\[
\mathscr{M}_A(Q_A(\mathbb{A}_AX)),M)
                     \iso\mathscr{N}_A(\mathbb{A}_AX,Z(M))
                     \iso\mathscr{M}_A(X,M),
\]
where the second isomorphism is a consequence of the universal
property of $\mathbb{A}_AX$.
\end{proof}

We will need to know the value of $\HAQ_*$ on sphere objects, and
here is the basic result required for this. This was proved in the
case of $A=S$ in lemma~3.6 and example~3.8 of Kuhn~\cite{Kuhn-TAQ}
using stabilisation, and the base change follows from the above
results.
For completeness, we present a more direct proof along similar
lines to that given by McCarthy and Minasian~\cite{RM&VM} in the
proof of theorem~6.1; unfortunately their argument appears to be
incorrect as stated (at one stage they incorrectly assume that~$M$
is an algebra).
\begin{proposition}\label{prop:Omega-PX/A}
Let $X$ be a cofibrant $A$-module, so that $\mathbb{P}_AX$ is a
cofibrant $A$-algebra. Then we have
\[
\Omega_A(\mathbb{P}_AX) \simeq \mathbb{P}_AX \wedge_A X.
\]
\end{proposition}
\begin{proof}
For every $M\in\mathscr{M}_{\mathbb{P}_A X}$ there is an adjunction
\[
\mathscr{C}_A / \mathbb{P}_A X (\mathbb{P}_A X,\mathbb{P}_A X\vee M)
 \iso \mathscr{M}_A / \mathbb{P}_A X(X,\mathbb{P}_A X\vee M),
\]
where $\mathscr{M}_A / \mathbb{P}_A X$ denotes the category of $A$-modules
over $\mathbb{P}_A X$. Because the forgetful functor
$\mathscr{C}_A / \mathbb{P}_A X \lra\mathscr{M}_A / \mathbb{P}_A X$
respects fibrations and acyclic fibrations, the adjunction passes
to homotopy categories. Now we have
\[
\mathscr{M}_A /\mathbb{P}_A X (\mathbb{P}_A X, M)
                                 \iso \mathscr{M}_A / X(X,X\vee M)
\]
and this passes to an adjunction between homotopy categories,
\[
\bar{h}\mathscr{M}_A /\mathbb{P}_A X (\mathbb{P}_A X, M)
                         \iso \bar{h}\mathscr{M}_A / X(X,X\vee M).
\]
Since in the homotopy category $X \vee M$ is the product of $X$
and $M$, we have
\[
\bar{h}\mathscr{M}_A / X(X,X\vee M) \iso \bar{h}\mathscr{M}_A(X,M).
\]
By using the free functor from $A$-modules to $\mathbb{P}_A X$-modules,
we obtain
\[
\bar{h}\mathscr{M}_A / X(X,X\vee M) \iso 
\bar{h} \mathscr{M}_{\mathbb{P}_A X}(\mathbb{P}_A X \wedge_A X, M).
\]
Thus we have shown that
\begin{align*}
\bar{h}\mathscr{M}_{\mathbb{P}_A X}(\Omega_A(\mathbb{P}_A X), M)
 &\iso \bar{h}\mathscr{C}_A / \mathbb{P}_A X(\mathbb{P}_A X ,\mathbb{P}_A X \vee M) \\
 &\iso \bar{h} \mathscr{M}_{\mathbb{P}_A X}(\mathbb{P}_A X \wedge_A X, M).
\end{align*}
Now using Yoneda's lemma, we have the desired equivalence.
\[
\Omega_A(\mathbb{P}_A X)\simeq \mathbb{P}_A X \wedge_A X.
\qedhere
\]
\end{proof}

On taking $B=\mathbb{P}_AX$ and $C=A$, the cofibration sequence
of~\cite{Basterra-TAQ}*{proposition~4.2} yields the cofibration
sequence of $A$-modules
\[
\Omega_A(\mathbb{P}_AX)\wedge_{\mathbb{P}_AX}A\lra\Omega_A(A)
                                  \lra\Omega_{\mathbb{P}_AX}(A)
\]
in which $\Omega_A(A)\simeq*$. Hence as $A$-modules,
\begin{equation}\label{eqn:Omega-A/PX}
\Omega_{\mathbb{P}_AX}(A)\simeq
\Sigma\Omega_A(\mathbb{P}_AX)\wedge_{\mathbb{P}_AX}A \simeq \Sigma X,
\end{equation}
where the second equivalence comes from Proposition~\ref{prop:Omega-PX/A}.
\begin{proposition}\label{prop:HAQ-PS^n/A}
For any $\mathbb{P}_AS^n$-module $M$ we have
\[
\TAQ_*(\mathbb{P}_AS^n,A;M)\iso M_{*-n},
\quad
\TAQ^*(\mathbb{P}_AS^n,A;M)\iso M^{*-n}.
\]
In particular,
\[
\HAQ_k(\mathbb{P}_AS^n,A)=\HAQ^k(\mathbb{P}_AS^n,A)=
\begin{cases}
\k& \text{if $k=n$}, \\
0& \text{otherwise}.
\end{cases}
\]
\end{proposition}
\begin{proof}
Taking $X=S^n$ in Proposition~\ref{prop:Omega-PX/A} we obtain
\begin{align*}
\TAQ_*(\mathbb{P}_AS^n,A;M)
&=\pi_*\Omega_{A}(\mathbb{P}_AS^n)\wedge_{\mathbb{P}_AS^n}M
 =\pi_*S^n\wedge M\iso M_{*-n}, \\
\TAQ^*(\mathbb{P}_AS^n,A;M)
&=\pi_{-*}F_{\mathbb{P}_AS^n}(\Omega_{A}(\mathbb{P}_AS^n),M)
 =\pi_{-*}F(S^n,M)\iso M^{*-n}.
\end{align*}
When $M=H\k$ with $\k=\pi_0A$ this gives
\[
\HAQ_k(\mathbb{P}_AS^n,A)=\HAQ^k(\mathbb{P}_AS^n,A)
=
\begin{cases}
\k& \text{\rm if $k=n$}, \\
0& \text{\rm otherwise},
\end{cases}
\]
as claimed.
\end{proof}
\begin{proposition}\label{prop:HAQ-A/PS^n}
We have
\[
\HAQ_k(A,\mathbb{P}_AS^n)=\HAQ^k(A,\mathbb{P}_AS^n)
=
\begin{cases}
\k& \text{\rm if $k=n+1$}, \\
0& \text{\rm otherwise}.
\end{cases}
\]
\end{proposition}
\begin{proof}
Taking $X=S^n$ in~\eqref{eqn:Omega-A/PX} we have
\begin{align*}
\Omega_{\mathbb{P}_AS^n}(A)\wedge_AH\k
&\simeq\Sigma\Omega_{A}(\mathbb{P}_AS^n)\wedge_{\mathbb{P}_AS^n}A\wedge_AH\k
 \simeq\Sigma S^n\wedge_AH\k, \\
F_{A}(\Omega_{\mathbb{P}_AS^n}(A),H\k)
&\simeq F_A(\Sigma\Omega_{A}(\mathbb{P}_AS^n)\wedge_{\mathbb{P}_AS^n}A,H\k)
 \simeq F_A(\Sigma S^n,H\k).
\end{align*}
Now using Proposition~\ref{prop:HAQ-PS^n/A}, the result is immediate.
\end{proof}

We say that the $S$-algebra $A$ is \emph{simply connected} if it
is connective and $\pi_0A=\Z$. Our next result is an analogue of
a standard result on connective spectra. Here we are considering
$\HAQ_*(A,S)=\HAQ_*(A,S;\Z)$.
\begin{proposition}\label{prop:HAQ-Whitehead}
Let $\phi\:A\lra B$ be a map of simply connected commutative $S$-algebras.
Then $\phi$ is an equivalence if and only if\/
$\phi_*\:\HAQ_*(A,S)\lra\HAQ_*(B,S)$ is an isomorphism. In particular,
the unit $S\lra B$ of such an $S$-algebra is a weak equivalence if
and only if $\HAQ_*(B,S)=0$.
\end{proposition}
\begin{proof}
Let $n\geq0$. Then $\phi$ is an $n$-equivalence if and only if the
mapping cone $\cone_\phi$ is $n$-connected. But on combining
Corollary~\ref{cor:Basterra8.2} with the long exact sequence
of~\eqref{eqn:HAQ_*}, we see that $\cone_\phi$ is $n$-connected if
and only if $\phi_*\:\HAQ_k(A,S)\lra\HAQ_k(B,S)$ is an isomorphism
for all $k\leq n$. Since this holds for all~$n$, the result follows.
\end{proof}

As a corollary we have an analogue of the Hurewicz isomorphism theorem.
\begin{corollary}\label{cor:HAQ-HurewiczThm}
Let $A$ be a commutative $S$-algebra whose unit $\eta\:S\lra A$ is an
$n$-equivalence. Then the Hurewicz homomorphism
$\theta\:\pi_{n+1}A\lra\HAQ_{n+1}(A,S)$ induces a monomorphism
\[
\theta'\:\pi_{n+1}A/\eta_*\pi_{n+1}S \lra \HAQ_{n+1}(A,S)
\]
which is an isomorphism if $\eta_*\:\pi_nS\lra\pi_nA$ is a monomorphism.
\end{corollary}

Here is a slightly different way to interpret the $\TAQ$ Hurewicz
homomorphism. Let $S^n_A$ denote the free $A$-module generated by
the sphere $S^n$. For a commutative $A$-algebra $B$ there for is
a Quillen adjunction
\[
\xymatrix{
\mathscr{M}_A(S^n_A,B)\ar@/^/@<-2ex>[r]
          & \mathscr{C}_A(\mathbb{P}_A S^n_A,B)\ar@/^/@<-2ex>[l]
}
\]
which gives rise to an isomorphism
\[
\pi_nB=\bar{h}\mathscr{M}_A(S^n_A,B)
                    \iso\bar{h}\mathscr{C}_A(\mathbb{P}_A S^n_A,B).
\]
So to each homotopy class $[f]\in\pi_nB$ we may assign the homotopy
class of $f'\in\mathscr{C}_A(\mathbb{P}_A S^n_A,B)$. Then for any
$B$ ring spectrum $E$, the induced homomorphism
\[
\Theta\:A_0\iso\pi_nS^n_A \lra E^A_nS^n_A\xrightarrow{\iso}
      \TAQ_n(\mathbb{P}_A S^n,A;E)\xrightarrow{f'_*}\TAQ_n(B,A;E)
\]
is related to the $\TAQ$ Hurewicz homomorphism by
\begin{equation}\label{eqn:Hurewicz}
\theta\:\pi_nB \lra \TAQ_n(B,A;E),
\quad
\theta([f]) = \Theta(1).
\end{equation}
A particularly important case occurs when $A=S$, and $E=H\pi_0B$.

\section{Topological Andr\'e-Quillen homology for cell $S$-algebras}
\label{sec:TAQhomology-cell}

We will apply the results of Section~\ref{sec:TAQhomology-recall} to
the case of a CW commutative $S$-algebra $R$ which is the colimit of
a sequence of cofibrations of cofibrant commutative $S$-algebras:
\[
S=R_{[0]} \xrightarrow{i_0} \cdots
\xrightarrow{i_{n-1}}R_{[n]}\xrightarrow{i_{n}}R_{[n+1]}\xrightarrow{i_{n+1}}
\cdots.
\]
We will also assume that only cells of degree greater than~$1$ are
attached, thus $R_{[1]}=R_{[0]}=S$ and $\pi_0R=\pi_0S$, so these
$S$-algebras are simply connected. The $(n+1)$-skeleton $R_{[n+1]}$
is obtained by attaching a wedge of $(n+1)$-cells to $R_{[n]}$ using
a map $k_n\:K_n\lra R_{[n]}$ from a wedge of $n$-spheres $K_n$ and
its extension to a map of $S$-algebras $\mathbb{P}_SK_n \lra R_{[n]}$.
To make this work properly, we need to use cofibrant replacement at
appropriate places.

So assume as an induction hypothesis inductively that $R_{[n]}$ is
defined and is a cofibrant $S$-algebra. Let $R'_{[n]}$ be a cofibrant
replacement for $R_{[n]}$ as an $\mathbb{P}_SK_n$-algebra, thus
there is an acyclic fibration $R'_{[n]}\twoheadrightarrow R_{[n]}$.

Now define $R_{[n+1]}$ and $R''_{[n+1]}$ to be the pushouts of the
diagrams
\begin{equation}\label{eqn:Pushout-R_(n+1)}
\xymatrix{
&{\mathbb{P}_SK_n}\ar[dl]\ar[dr]& \\
{\mathbb{P}_S\mathrm{C}K_n}&&R_{[n]}
}
\qquad
\xymatrix{
&{\mathbb{P}_SK_n}\ar[dl]\ar[dr]& \\
{\mathbb{P}_S\mathrm{C}K_n}&&R'_{[n]}
}
\end{equation}
which are given by
\begin{equation}\label{eqn:R^(n+1)-formula}
R_{[n+1]}=\mathbb{P}_S\mathrm{C}K_n\wedge_{\mathbb{P}_SK_n}R_{[n]},
\quad
R''_{[n+1]}=\mathbb{P}_S\mathrm{C}K_n\wedge_{\mathbb{P}_SK_n}R'_{[n]}.
\end{equation}
Now consider the diagram of commutative $\mathbb{P}_SK_n$-algebras
(and hence of $S$-algebras)
\[
\xymatrix{
&{\mathbb{P}_SK_n}\ar@{ >->}[dl]\ar@{ >->}[dr]& & \\
{\mathbb{P}_S\mathrm{C}K_n}\ar[dr]&&R'_{[n]}\ar@{ >->}[dl]\ar@{->>}[dr]^{\sim} &\\
&R''_{[n+1]}\ar[dr]_{\sim}& &R_{[n]}\ar@{ >->}[dl] \\
&&R_{[n+1]} &
}
\]
in which the upper and composite parallelograms are pushout
diagrams. Since the top left hand arrow is a cofibration,
the parallel ones are also cofibrations as they are pushouts.
Also, by~\cite{EKMM}*{proposition~VII.7.4} the functor
${\mathbb{P}_S\mathrm{C}K_n}\wedge_{\mathbb{P}_SK_n}(-)$
preserves weak equivalences between cofibrant $S$-algebras,
therefore the lower left arrow is also a weak equivalence.

Notice that for each~$n$, there is a lifting diagram of the
form
\[
\xymatrix{
S\ar[r]\ar@{ >->}[d] & R'_{[n]}\ar@{->>}[d]^{\sim} \\
R_{[n-1]}\ar[r]\ar@{.>}[ur]              & R_{[n]}
}
\]
so it does no harm to assume that at each stage we have replaced
$R_{[n]}$ by $R'_{[n]}$ in what follows.

Now by~\cite{Basterra-TAQ}*{proposition~4.6}, for a cofibrant
$S$-algebra $A$ and cofibrant $A$-algebras $A\lra B$ and $A\lra C$,
we have
\begin{equation}\label{eqn:Omega-Change}
\Omega_{C}(B\wedge_A C)\simeq\Omega_A(B)\wedge_A C.
\end{equation}
For $n\geq1$ this gives
\[
\Omega_{R_{[n]}}(\mathbb{P}_S\mathrm{C}K_n\wedge_{\mathbb{P}_SK_n}R_{[n]})
\simeq
\Omega_{\mathbb{P}_SK_n}(\mathbb{P}_S\mathrm{C}K_n)\wedge_{\mathbb{P}_SK_n}R_{[n]}
\]
and hence there is a long exact sequence derived from~\eqref{eqn:TAQ-LES_*},
\begin{multline*}
\cdots\lra\HAQ_k(R_{[n]},S)\lra\HAQ_k(R_{[n+1]},S)
\lra\HAQ_k(R_{[n+1]},R_{[n]}) \\
\lra\HAQ_{k-1}(R_{[n]},S)\lra\cdots
\end{multline*}
which by~\eqref{eqn:Omega-Change} becomes
\begin{multline*}
\cdots\lra\HAQ_k(R_{[n]},S)\lra\HAQ_k(R_{[n+1]},S)
\lra\HAQ_k(\mathbb{P}_S\mathrm{C}K_n,\mathbb{P}_SK_n) \\
\lra\HAQ_{k-1}(R_{[n]},S)\lra\cdots
\end{multline*}
in which there is an equivalence of $\mathbb{P}_SK_n$-algebras
\[
\mathbb{P}_S\mathrm{C}K_n\simeq\mathbb{P}_S*=S.
\]
Hence we obtain the following long exact sequence
\begin{multline}\label{eqn:CWR-LES}
\cdots\lra\HAQ_k(R_{[n]},S)\lra\HAQ_k(R_{[n+1]},S)
\lra\HAQ_k(S,\mathbb{P}_SK_n)  \\
\lra\HAQ_{k-1}(R_{[n]},S)\lra\cdots.
\end{multline}

Using Proposition~\ref{prop:HAQ-A/PS^n}, we can now give an
estimate for the size of $\HAQ_*(R,S)$ when~$R$ is a finite
dimensional CW~commutative $S$-algebra.
\begin{proposition}\label{prop:HAQ-finCW}
Let $R$ be a CW commutative $S$-algebra with cells only in
degrees at most~$n$. Then for $k>n$, $\HAQ_k(R,S)=0$.
\end{proposition}
\begin{corollary}\label{cor:HAQ-finCW}
If~$R$ has only finitely many cells, then
\[
\sum_{k=0}^n\rank\HAQ_k(R,S)\leq\text{\rm number of cells}.
\]
\end{corollary}

In the category of $S$-modules there are three cofibration
sequences that will concern us. We have the two cofibration
sequences
\[
K_n\xrightarrow{k_n}R_{[n]}\lra\cone_{k_n},
\quad
R_{[n]}\xrightarrow{i_n}R_{[n+1]}\lra\cone_{i_n}.
\]
{}From the proof of~\cite{Basterra-TAQ}*{lemma~8.2}, there
is a homotopy commutative diagram
\[
\xymatrix{
R_{[n]}\ar[r]^{i_n}\ar[d] & R_{[n+1]}\ar[r]\ar[d] & \cone_{i_n}\ar[d]^{\tau_n} \\
R_{[n]}\ar[r]^{i_n} & R_{[n+1]}\ar[r]^{\delta_n\ph{-}} & \Omega_{R_{[n]}}(R_{[n+1]})
}
\]
where $\delta_n$ denotes the universal derivation; we claim
this extends to a homotopy commutative diagram of the following
form.
\begin{equation}\label{eqn:SettingUp-1}
\xymatrix{
R_{[n]}\ar[r]\ar[d]^{=} & \cone_{k_n}\ar[r]\ar[d] & \Sigma K_n\ar[d]^{h_n} \\
R_{[n]}\ar[r]^{i_n}\ar[d]^{=} &R_{[n+1]}\ar[r]\ar[d]^{=} &\cone_{i_n}\ar[d]^{\tau_n} \\
R_{[n]}\ar[r]^{i_n} &R_{[n+1]}\ar[r]^{\delta_n\ph{-}} &\Omega_{R_{[n+1]}/R_{[n]}}
}
\end{equation}
Recalling that $R_{[n+1]}$ is a pushout for the diagram of commutative
$S$-algebras~\eqref{eqn:Pushout-R_(n+1)}, the map $\cone_{k_n}\lra R_{[n+1]}$
exists since $\cone_{k_n}$ is defined as a pushout for the diagram of
$S$-modules
\[
\xymatrix{
&K_n\ar[dl]_{k_n}\ar[dr]& \\
R_{[n]}&&{\mathrm{C}K_n}
}
\]
and the evident composition
\[
\mathrm{C}K_n\lra\mathbb{P}_S\mathrm{C}K_n\lra R_{[n+1]}
\]
gives rise to a commutative diagram of $S$-modules of the
following form.
\[
\xymatrix{
&K_n\ar[dl]_{k_n}\ar[dr]& \\
R_{[n]}\ar[dr]_{i_n}&&{\mathrm{C}K_n}\ar[dl] \\
&R_{[n+1]}&
}
\]
Using~\eqref{eqn:Omega-Change}, we see that there are
equivalences of modules over
$R_{[n+1]} = R_{[n]}\wedge_{\mathbb{P}_SK_n}\mathbb{P}_SCK_n$:
\begin{equation}\label{eqn:OmegaR[n+1]/R[n]}
\Omega_{R_{[n]}}(R_{[n+1]})
\simeq
R_{[n]}\wedge_{\mathbb{P}_SK_n}\Omega_{\mathbb{P}_SK_n}(S)
\simeq
R_{[n]}\wedge_{\mathbb{P}_SK_n}S\wedge\Sigma K_n.
\end{equation}
Now smashing over $R_{[n+1]}$ with $H\Z$ we obtain
\begin{equation}\label{eqn:HZOmegaR[n+1]/R[n]}
H\Z\wedge_{R_{[n+1]}}\Omega_{R_{[n]}}(R_{[n+1]})\simeq H\Z\wedge\Sigma K_n.
\end{equation}
On smashing $\tau_n\o h_n$ with $H\Z$ we obtain a map
\[
H\Z\wedge\Sigma K_n\xrightarrow{\;\id\wedge\tau_n\o h_n\;}
H\Z\wedge\Omega_{R_{[n]}}(R_{[n+1]})
\]
and following this with the natural map
\begin{equation}\label{eqn:fn-Equivalence}
H\Z\wedge\Omega_{R_{[n]}}(R_{[n+1]})
   \lra H\Z\wedge_{R_{[n+1]}}\Omega_{R_{[n]}}(R_{[n+1]})
   \simeq H\Z\wedge\Sigma K_n
\end{equation}
yields a self map $f_n\:H\Z\wedge\Sigma K_n\lra H\Z\wedge\Sigma K_n$.
Since $K_n$ is a wedge of $n$-spheres, $H\Z\wedge K_n$ is a wedge
of copies of $H\Z$. In fact, the map of~\eqref{eqn:fn-Equivalence}
induces an isomorphism on $\pi_{n+1}(-)$.
\begin{lemma}\label{lem:fn-Equivalence}
The map $f_n\:H\Z\wedge\Sigma K_n\lra H\Z\wedge\Sigma K_n$ is a
weak equivalence. Equivalently, the following maps are isomorphisms:
\[
\pi_{n+1}\Sigma K_n\xrightarrow{(h_n)_*}\pi_{n+1}\cone_{i_n},
\quad
\pi_{n+1}\Sigma K_n\xrightarrow{(\tau_n\o h_n)_*}
\pi_{n+1}\Omega_{R_{[n]}}(R_{[n+1]}).
\]
\end{lemma}
\begin{proof}
The pairs $(\cone_{k_n},R_{[n]})$ and $(R_{[n+1]},R_{[n]})$ occurring
in~\eqref{eqn:SettingUp-1} are relative cell complexes which have the
same cells in degrees up to $2n+1$. The cells in degree $n+1$ correspond
to those on $\Sigma K_n$ and therefore
$(h_n)_*\:\pi_{n+1}\Sigma K_n\lra\pi_{n+1}\cone_{i_n}$ is an isomorphism.
For a discussion of cellular structures in this context,
see~\cite{EKMM}*{VII~3,~X 2}.

It now follows from the Hurewicz isomorphism theorem that $f_n$ induces
an isomorphism on $\pi_{n+1}(H\Z\wedge \Sigma K_n)$ which agrees with
$H_{n+1}(\Sigma K_n)$.
\end{proof}

Applying homotopy to the diagram of~\eqref{eqn:SettingUp-1}, we obtain
a diagram of abelian groups, a part of which is
\begin{equation}\label{eqn:SettingUp-2}
\xymatrix{
\cdots \ar[r] & \pi_{n+1}R_{[n]}\ar[r]\ar[d]_{=}
               & \pi_{n+1}\cone_{k_n} \ar[r]\ar[d]
               &\pi_{n+1}\Sigma K_n \ar[r]\ar[d]_{\iso}^{(h_n)_*}
               & \pi_nR_{[n]}\ar[r]\ar[d]_{=}
                                       & \cdots \\
\cdots\ar[r] &\pi_{n+1}R_{[n]}\ar[r]^{(i_n)_*}\ar[d]_{=}
             &\pi_{n+1}R_{[n+1]}\ar[r]\ar[d]_{=}
             &\pi_{n+1}\cone_{i_n}\ar[r]\ar[d]_{\iso}^{(\tau_n)_*}
             & \pi_{n}R_{[n]}\ar[r]\ar[d]_{=}
                                       & \cdots \\
\cdots \ar[r] & \pi_{n+1}R_{[n]} \ar[r]^{(i_n)_*}
             & \pi_{n+1}R_{[n+1]}\ar[r]^{(\delta_n)_*\ph{--}}
             & \pi_{n+1}\Omega_{R_{[n]}}(R_{[n+1]})\ar[r]
             & \pi_nR_{[n]}\ar[r] & \cdots
}
\end{equation}
whose top two rows are exact. In the portion shown, the bottom row
is also exact because $(\tau_n)_*$ is an isomorphism on $\pi_{n+1}(-)$.

Using the definition in terms of~\eqref{eqn:Hurewicz-UnivDeriv-E},
the Hurewicz homomorphism
\[
\theta_{n+1}\:\pi_{n+1}R_{[n+1]}\lra\HAQ_{n+1}(R_{[n+1]},R_{[n]})
\]
is induced from the composition
\[
R_{[n+1]}\xrightarrow{\delta_{(R_{[n+1]},R_{[n]})}} \Omega_{R_{[n]}}(R_{[n+1]})
   \lra H\Z\wedge_{R_{[n]}}\Omega_{R_{[n]}}(R_{[n+1]})
   \lra H\Z\wedge_{R_{[n+1]}}\Omega_{R_{[n]}}(R_{[n+1]}),
\]
and using~\eqref{eqn:Hurewicz-diagram} it extends to a diagram
\begin{equation}\label{eqn:SettingUp-3}
\xymatrix{
\pi_{n+1}R_{[n+1]}\ar[r]\ar[d]^{\theta_{n+1}}
 & \pi_{n+1}\Omega_{R_{[n]}}(R_{[n+1]})\ar[r]\ar[d]^{\iso}
 & \pi_nR_{[n]}\ar[d]^{\theta_n} \\
\HAQ_{n+1}(R_{[n+1]},S)\ar[r]
& \HAQ_{n+1}(R_{[n+1]},R_{[n]})\ar[r]
& \HAQ_n(R_{[n]},S)
}
\end{equation}
in which the bottom row is a portion of the usual long exact
sequence~\eqref{eqn:TAQ-LES} for $A=S$, $B=R_{[n]}$ and $C=R_{[n+1]}$.
Furthermore, these diagrams are compatible for varying~$n$.

Using the evident natural transformation $\HAQ_n(-)\lra\HAQ_n(-;\F_p)$,
we can map the bottom row of~\eqref{eqn:SettingUp-3} into the exact
sequence
\[
0\ra\HAQ_{n+1}(R_{[n+1]},S;\F_p) \lra
              \HAQ_{n+1}(R_{[n+1]},R_{[n]};\F_p)
                             \lra\HAQ_n(R_{[n]},S;\F_p)
\]
to obtain the commutative diagram
\begin{equation}\label{eqn:SettingUp-5}
\xymatrix{
\pi_{n+1}R_{[n+1]}\ar[r]\ar[d]^{\bar\theta_{n+1}}
 & \pi_{n+1}\Sigma K_n\ar[r]\ar[d]^{\text{epi}}
 & \pi_nR_{[n]}\ar[d]^{\bar\theta_n}  \\
0\ra\HAQ_{n+1}(R_{[n+1]},S;\F_p)\ar[r]
 & \HAQ_{n+1}(R_{[n+1]},R_{[n]};\F_p)\ar[r]
 & \HAQ_n(R_{[n]},S;\F_p)
}
\end{equation}
in which the rows are exact and the middle vertical arrow
is an epimorphism.
\begin{remark}\label{rem:S->R_0}
All of the above works as well if we replace $S$ by the
$p$-local sphere for some prime~$p$, or more generally by
a connective commutative $S$-algebra $A$ with $\pi_0A$ a
localisation of $\pi_0S$.
\end{remark}

\section{Minimal atomic $p$-local commutative $S$-algebras}
\label{sec:MinAtomic}

{}From now on we fix a prime~$p$. We work with $p$-local spectra,
suppressing indications of localization from the notation when
convenient. Thus $S$ stands for $S_{(p)}$ and so on. When discussing
cell structures we always mean these to be taken $p$-locally, for
example in the category of $S$-modules or $S$-algebras.

In~\cite{AJB-JPM}, the following notion was introduced. A $p$-local
CW~complex $Y$ is \emph{minimal} if its cellular chain complex with
$\F_p$ coefficients has trivial boundaries, so
\[
\mathrm{H}_*(Y;\F_p)=\mathrm{C}_*(Y)\otimes\F_p.
\]
An alternative formulation in terms of the skeletal inclusion maps
$Y_n\lra Y_{n+1}\lra Y$ is that for each $n$ the induced epimorphism
\[
\mathrm{H}_n(Y_n;\F_p)\lra\mathrm{H}_n(Y_{n+1};\F_p)
\]
is actually an isomorphism and so
\[
\mathrm{H}_n(Y_n;\F_p)\xrightarrow{\iso}\mathrm{H}_n(Y_{n+1};\F_p)
\xrightarrow{\iso}\mathrm{H}_n(Y;\F_p).
\]
In~\cite{AJB-JPM}*{theorem~3.3} it was shown that every $p$-local
CW~complex of finite-type is equivalent to a minimal one, so such
minimal complexes exist in abundance.

Continuing to take the point of view that $\HAQ_*$ is a good substitute
for ordinary homology when considering commutative $S$-algebras, let us
consider the analogous notion in this multiplicative situation. We begin
with a suitable definition of minimality in this context.
\begin{definition}\label{defn:MinimalCW}
Let $R$ be a $p$-local CW commutative $S$-algebra with $n$-skeleton
$R_{[n]}$. Then $R$ is \emph{minimal} if for each~$n$ the inclusion
maps $R_{[n]}\lra R_{[n+1]}\lra R$ induce isomorphisms
\[
\HAQ_n(R_{[n]},S;\F_p)\xrightarrow{\iso}\HAQ_n(R_{[n+1]},S;\F_p)
\]
or equivalently
\[
\HAQ_n(R_{[n]},S;\F_p)\xrightarrow{\iso}\HAQ_n(R,S;\F_p).
\]
\end{definition}

Proposition~\ref{prop:HAQ-finCW} implies that the homomorphisms
here are both epimorphisms, as is true for their analogues in
ordinary homology.
\begin{theorem}\label{thm:MinimalCW-Existence}
Let $R$ be $p$-local CW commutative $S$-algebra, with finitely
many $p$-local cells in each degree. Then there is a minimal
$p$-local CW~commutative $S$-algebra $R'$ and an equivalence
of commutative $S$-algebras $R'\lra R$.
\end{theorem}
\begin{proof}
The details follow from the proof of~\cite{AJB-JPM}*{theorem~3.3},
replacing ordinary homology with $\HAQ_*(-)$, \emph{mutatis mutandis}.
\end{proof}

This allows us to revisit the results of~\cite{AJB-JPM} in the context
of simply connected $p$-local CW~commutative $S$-algebras. The notion
of a \emph{nuclear algebra} appears in~\cite{Hu-Kriz-May}*{definition~2.7},
while that of an \emph{atomic algebra} appears
in~\cite{Hu-Kriz-May}*{definition~2.8}. We also use notions
from~\cite{AJB-JPM}*{definition~1.1}. Further details and related
results appeared in the second author's PhD thesis~\cite{HG:PhD}.

First we have an analogue of~\cite{AJB-JPM}*{theorem~3.4}. We remark
that it was not pointed out explicitly in~\cite{AJB-JPM} that a nuclear
complex is always minimal so our next result has a direct analogue in
the context of the earlier work even though it did not appear in the
published version.
\begin{theorem}\label{thm:Nuc<->min&nomodphtpy}
Let $R$ be a simply connected $p$-local CW commutative $S$-algebra.
Then~$R$ is nuclear if and only if it is minimal and the Hurewicz
homomorphism $\theta\:\pi_nR\lra\HAQ_n(R,S;\F_p)$ is trivial for
all $n>0$.
\end{theorem}
\begin{proof}
Suppose that $R$ is nuclear in the discussion of
Section~\ref{sec:TAQhomology-cell}. Then as there is a factorisation
\[
\xymatrix{
\pi_{n+1}R_{[n+1]}\ar[rr]\ar[dr]&&\pi_{n+1}\Sigma K_n \\
&p\,\pi_{n+1}\Sigma K_n\ar[ur]&
}
\]
we have $\bar\theta_{n+1}=0$ for $n\geq0$, and also
\[
\HAQ_0(R_{[0]},S;\F_p)=\HAQ_0(S,S;\F_p)=0.
\]
Similarly, the image of the boundary map
$\HAQ_{n+1}(R_{[n+1]},R_{[n]};\F_p)\lra\HAQ_n(R_{[n]},S;\F_p)$ is
contained in $\im\bar\theta_n=0$, hence it is trivial. This shows
that $R$ is minimal in the sense of Definition~\ref{defn:MinimalCW},
and has no mod~$p$ detectable homotopy in the sense that for $n>0$,
the composition
\[
\bar\theta_n\:\pi_nR_{[n]}\lra\HAQ_n(R_{[n]},S;\F_p)\lra\HAQ_n(R,S;\F_p)
\]
is trivial.

By a similar argument to the proof of~\cite{AJB-JPM}*{theorem~3.4},
the converse also holds, \ie, if $R$ is minimal and has no mod~$p$
detectable homotopy then it is nuclear.
\end{proof}

We claim that~\cite{Hu-Kriz-May}*{conjecture~2.9} and the analogue
of~\cite{AJB-JPM}*{proposition~2.5} are consequences of our next
result.
\begin{theorem}\label{thm:Nuc->Atomic}
A nuclear simply connected $p$-local CW commutative $S$-algebra is
minimal atomic, hence a core of such an algebra is an equivalence.
\end{theorem}
\begin{proof}
This follows the analogous proof of~\cite{AJB-JPM}.
\end{proof}

We also note the following useful result.
\begin{proposition}\label{prop:minatomS-mod->A-alg}
Let $R$ be minimal atomic as an $S$-module. Then it is minimal
atomic as an $S$-algebra.
\end{proposition}
\begin{proof}
This follows easily from the fact that for $S$-modules, being
minimal atomic is equivalent to being irreducible.

As an alternative, notice that if $R$ is minimal atomic as an
$S$-module, then the ordinary homology Hurewicz homomorphism
$\pi_nR\lra H_n(R;\F_p)$ is trivial for $n>0$, so by
Proposition~\ref{prop:Hurewicz-E} the mod~$p$ $\HAQ$ Hurewicz
homomorphism $\pi_nR\lra\HAQ_n(R,S;\F_p)$ is trivial, whence
$R$ is a minimal atomic $S$-algebra.
\end{proof}

\section{The $\TAQ$ Hurewicz homomorphism for Thom spectra}
\label{sec:Hurewicz-Thomspec}

In order to calculate with Thom spectra arising from infinite
loop maps into $BSF$, we need some information on the relevant
universal derivations. The next two results are implicit in the
proof of~\cite{Basterra-Mandell}*{theorem~6.1}, but unfortunately
they are not stated explicitly and we are grateful to Mike Mandell
and Maria Basterra for clarifying this material which is due to
them.

Let $X$ be an infinite loop space and let $\underline{X}$ be the
associated spectrum viewed as an $S$-module.
Following~\cite{Basterra-Mandell}*{section~6}, consider the
augmented $S$-algebra
\[
\Sigma^\infty_{S+}X = S\wedge_{\mathcal{L}}\Sigma^\infty X_+.
\]
There is a canonical (evaluation) morphism of $S$-modules
$\sigma\:\Sigma^\infty_{S+}X\lra\underline{X}$. Now taking $A=B=S$
and using the isomorphism
$\Sigma^\infty_{S+}X\wedge_S S\iso\Sigma^\infty_{S+}X$, the bijection
of~\cite{Basterra-TAQ}*{proposition~3.2} gives
\begin{equation}\label{eqn:MBprop3.2}
\bar{h}\mathscr{C}_S/S(\Sigma^\infty_{S+}X,S\vee M)
\iso
\bar{h}\mathscr{M}_S(\mathrm{L}Q_S\mathrm{R}I_S(\Sigma^\infty_{S+}X),M).
\end{equation}
When $M=\mathrm{L}Q_S\mathrm{R}I_S(\Sigma^\infty_{S+}X)$, the identity
morphism on the right hand side corresponds to a morphism
\[
\Sigma^\infty_{S+}X\lra
          S\vee\mathrm{L}Q_S\mathrm{R}I_S(\Sigma^\infty_{S+}X)
\]
which projects to a universal derivation
\[
\delta_X\:\Sigma^\infty_{S+}X\lra
                     \mathrm{L}Q_S\mathrm{R}I_S(\Sigma^\infty_{S+}X).
\]
\begin{proposition}\label{prop:UnivDer-Suspspec}
In the homotopy category $\bar{h}\mathscr{M}_S$ there is an isomorphism
\[
\theta\:\underline{X}\xrightarrow{\iso}
   \mathrm{L}Q_S\mathrm{R}I_S(\Sigma^\infty_{S+}X)
\]
and a commutative diagram
\[
\xymatrix{
&\Sigma^\infty_{S+}X\ar[dl]_{\sigma}\ar[dr]^{\delta_X}& \\
\underline{X}\ar[rr]^{\theta\ph{ABCDE}}
&&
\mathrm{L}Q_S\mathrm{R}I_S(\Sigma^\infty_{S+}X).
}
\]
Hence $\sigma$ realises the universal derivation.
\end{proposition}

Note that the diagonal on $X$ induces a morphism of commutative
$\Sigma^\infty_{S+}X$-algebras
\[
\Sigma^\infty_{S+}X\wedge_S \Sigma^\infty_{S+}X
\xrightarrow{\mathrm{id}\wedge\mathrm{diag}}
\Sigma^\infty_{S+}X\wedge_S \Sigma^\infty_{S+}X\wedge_S \Sigma^\infty_{S+}X
\xrightarrow{\mathrm{mult}\wedge\mathrm{id}}
\Sigma^\infty_{S+}X\wedge_S \Sigma^\infty_{S+}X
\]
which is a weak equivalence. This has the effect of interchanging the
multiplication map with the augmentation onto the second factor, giving
rise to a composition of isomorphisms in
$\bar{h}\mathscr{M}_{\Sigma^\infty_{S+}X}$, namely
\begin{equation}\label{eqn:Switch}
\Omega_{S}(\Sigma^\infty_{S+}X) \xrightarrow{\iso}
\Sigma^\infty_{S+}X\wedge_S\mathrm{L}Q_S\mathrm{R}I_S(\Sigma^\infty_{S+}X)
\xrightarrow{\iso}
\Sigma^\infty_{S+}X\wedge_S\underline{X}.
\end{equation}

For an infinite loop map $f\:X\lra BF$, the resulting Thom spectrum~$Mf$
(viewed as a commutative $S$-algebra) has a Thom diagonal
$\Delta\:Mf\lra Mf\wedge_S\Sigma^\infty_{S+}X$, so that the composition
\[
Mf\wedge_S Mf\xrightarrow{\mathrm{id}\wedge\Delta}
Mf\wedge_S Mf\wedge_S\Sigma^\infty_{S+}X
    \xrightarrow{\mathrm{mult}\wedge\mathrm{id}}
                      Mf\wedge_S\Sigma^\infty_{S+}X
\]
is a weak equivalence of $Mf$-algebras. Furthermore, the map
$\sigma$ induces a map of $S$-algebras
\[
\delta_f\:Mf\xrightarrow{\Delta}Mf\wedge_S\Sigma^\infty_{S+}X
\xrightarrow{\mathrm{id}\wedge\sigma} Mf\wedge_S\underline{X}
\]
and then we have the following result in which we use the universal
derivation $\delta_{(Mf,S)}$ of~\eqref{eqn:Omega-UnivDer}.
\begin{proposition}\label{prop:UnivDer-Thomspec}
In the homotopy category $\bar{h}\mathscr{M}_{Mf}$ there is an
isomorphism
\[
\Theta\:Mf\wedge_S\underline{X}\xrightarrow{\ph{-}\iso\ph{-}}\Omega_{S}(Mf)
\]
and a commutative diagram
\[
\xymatrix{
&Mf\ar[dl]_{\delta_f}\ar[dr]^{\delta_{(Mf,S)}}& \\
Mf\wedge_S\underline{X}\ar[rr]^{\Theta}
&&
\Omega_{S}(Mf),
}
\]
and hence $\delta_f$ realises the universal derivation.
\end{proposition}

Given this result, we can now describe how to calculate the $\HAQ$
Hurewicz homomorphism for a Thom spectrum $Mf$ induced by an infinite
loop map $f\:X\lra BSF$. Recalling the definition using the map
of~\eqref{eqn:Hurewicz-UnivDeriv-E}, we can take $\TAQ$ with
coefficient spectrum the integral Eilenberg-Mac~Lane spectrum $H=H\Z$
(or indeed the Eilenberg-Mac~Lane spectrum of any commutative ring)
and so obtain homomorphisms
\[
\theta\:\pi_n Mf \lra \HAQ_n(Mf,S) = \TAQ_n(Mf,S;H).
\]
Then the map of~\eqref{eqn:Hurewicz-UnivDeriv-E} factors as in the
commutative diagram obtained from~\eqref{eqn:theta-defndiagram},
\begin{equation}\label{eqn:theta-Thom}
\xymatrix{
Mf \ar[r]^{\ph{A}\delta_{(Mf,S)}\ph{ABC}}\ar[d]& \Omega_{S}(Mf) \ar[r]\ar[d]
& H\wedge\Omega_{S}(Mf)\ar[r]\ar[d] & H\wedge_{Mf}\Omega_{S}(Mf) \\
H\wedge Mf \ar[r]^{\mathrm{id}\wedge\delta_{(Mf,S)\ph{-}}\ph{a}}
& H\wedge \Omega_{S}(Mf) \ar[r]& H\wedge H\wedge\Omega_{S}(Mf)\ar[r]
& H\wedge H\wedge_{Mf}\Omega_{S}(Mf)\ar[u]^{\mathrm{mult}\wedge\mathrm{id}}
}
\end{equation}
in which the undecorated smash products are taken over~$S$ and the
downward pointing arrows are obtained by smashing with the unit
$S\lra H$. Thus $\theta$ factors through the usual Hurewicz homomorphism:
\[
\theta\:\pi_nMf \lra H_n(Mf) \lra \HAQ_n(Mf,S).
\]
Using Proposition~\ref{prop:UnivDer-Thomspec} we see that $\theta$ is
equivalent to the homomorphism
\[
\theta'\:\pi_nMf \lra H_n(\underline{X}) = \pi_n(H\wedge\underline{X})
\]
induced by
\[
Mf\xrightarrow{\Delta} Mf\wedge \Sigma^\infty X_+
\xrightarrow{U\wedge\mathrm{id}} H\wedge\Sigma^\infty X_+
\xrightarrow{\mathrm{id}\wedge\sigma} H\wedge\underline{X},
\]
where $U\in H^0(Mf;R)$ is the orientation class. On smashing with $H$,
we obtain
\[
\xymatrix{
H\wedge Mf\ar[r]\ar@{-->}@/_3.7pc/[rrrd]\ar@{.>}@/_1pc/[rrd] &
H\wedge Mf\wedge \Sigma^\infty X_+\ar[r] &
H\wedge H\wedge\Sigma^\infty X_+\ar[r]\ar[d] &
H\wedge H\wedge H\wedge\underline{X}\ar[d] \\
&&H\wedge\Sigma^\infty X_+ \ar[r]&
H\wedge\underline{X} \\
&&&
}
\]
where the dashed arrow induces $\theta''\:H_n(Mf)\lra\HAQ_n(Mf,S)$
and the dotted arrow induces the Thom isomorphism $H_n(Mf)\lra H_n(X)$.
Thus $\theta''$ is the Thom isomorphism composed with
$\sigma_*\:H_n(X)\lra H_n(\underline{X})$.

Notice that $\sigma_*$ factors through the homology suspension
\[
H_n(X)\xrightarrow{\iso}H_{n+1}(\Sigma X)\lra H_{n+1}(BX),
\]
where $BX$ is the delooping of $X$. Hence $\sigma_*$, and therefore
$\theta''$, annihilates decomposables in the rings on which they
are defined.

\section{Some examples of Thom spectra}\label{sec:Examp-Thomspec}

\noindent
\textbf{The Thom spectrum $MU$:}
By Proposition~\ref{prop:UnivDer-Thomspec}, we know that
$\Omega_{S}(MU)=MU\wedge\Sigma^2 ku$ since $BU$ is the zeroth
space in the $1$-connected cover of $ku$ which is $\Sigma^2 ku$.
Thus for any commutative ring $\F$ we have
\[
\HAQ_*(MU,S;\F) = H_*(\Sigma^2 ku;\F) = H_{*-2}(ku;\F).
\]
For the prime $2$,
\[
\HAQ_*(MU,S;\F_2) = H_{*-2}(ku;\F_2),
\]
where the right hand side is given by
\[
H_*(ku;\F_2) = \F_2[\zeta_1^4,\zeta_2^2,\zeta_3,\zeta_4,\ldots]
\subseteq\mathcal{A}(2)_*=\F_2[\zeta_1,\zeta_2,\zeta_3,\zeta_4,\ldots],
\]
and is a sub Hopf algebra of the mod~$2$ dual Steenrod algebra
with the inclusion induced from the natural map $ku\lra H\F_2$.
The generators $\zeta_r$ are the conjugates of the generators
$\xi_r$ coming from $H_*(B\Z/2;\F_2)$.

For an odd prime $p$, there is decomposition
\[
ku_{(p)} \simeq \bigvee_{0\leq r\leq p-2} \Sigma^{2r}\ell,
\]
where $\ell$ is the Adams summand for which
\[
H_*(\ell;\F_p) =
\F_p[\zeta_1,\zeta_2,\zeta_3,\ldots]\otimes
                            \Lambda(\bar\tau_2,\bar\tau_3,\ldots)
\subseteq\mathcal{A}(p)_*=\F_p[\zeta_1,\zeta_2,\zeta_3,\ldots]\otimes
                            \Lambda(\bar\tau_0,\bar\tau_1,\ldots),
\]
where the embedding into the dual mod~$p$ dual Steenrod algebra
is induced by a map $\ell\lra H\F_p$. The generators $\zeta_r$
and $\bar\tau_r$ are the conjugates of the generators $\xi_r$
and $\tau_r$ coming from $H_*(B\Z/p;\F_p)$.

Recall that for a prime~$p$, the image of the mod~$p$ Hurewicz
homomorphism $h\:\pi_*MU \lra H_*(MU;\F_p)$ is
\[
h\pi_*MU=\F_p[x_r:\text{$r+1$ is not a power of $p$}]
                                        \subseteq H_*(MU;\F_p),
\]
where, in the notation of Adams~\cite{JFA:Chicago},
\[
x_r \equiv b_r \pmod{\mathrm{decomposables}}.
\]
Of course, $x_r$ is the image of an element of $\pi_{2r}MU_{(p)}$.
The Thom isomorphism is
\[
\Phi\:H_*(MU;\F_p) \xrightarrow{\iso} H_*(BU;\F_p);
\quad
b_n \longleftrightarrow \beta_n,
\]
so to determine the Hurewicz homomorphism
$\theta\:\pi_*MU \lra \HAQ_*(MU,S;\F_p)$ we need to calculate
the image of $\sigma_*\: H_*(BU;\F_p) \lra H_{*-2}(ku;\F_p)$.
Since decomposables are killed by $\sigma_*$, it suffices to
know its effect on the generators $\beta_n$.

When $p=2$, recall that the map $ku\lra H\F_2$ induces a
monomorphism in $H_*(-;\F_2)$ and the canonical maps give
a composition
\[
\Sigma^\infty \CP^\infty \lra \Sigma^\infty BU
                         \lra \Sigma^2ku \lra \Sigma^2H\F_2
\]
which corresponds to the natural cohomology generator
$x\in H^2(\CP^\infty;\F_2)$ and the homology generator $\beta_n$
dual to $x^n$ (which maps to $\beta_n\in H_{2n}(BU;\F_p)$) maps
to the coefficient of $t^n$ in the power series
\[
\xi(t) = \sum_{s\geq0} \xi_s^2 t^{2^s}.
\]
This shows that
\begin{equation}\label{eqn:sigma*}
\sigma_*(\beta_n)
\begin{cases}
\neq 0 & \text{if $n$ is a power of $2$}, \\
= 0 & \text{otherwise}.
\end{cases}
\end{equation}
Hence we see that
\begin{equation}\label{eqn:theta}
\theta\pi_{2n}MU
\begin{cases}
\neq 0 & \text{if $n$ is a power of $2$}, \\
= 0 & \text{otherwise}.
\end{cases}
\end{equation}
More precisely, this shows that the polynomial generator of $\pi_*MU$
in degree~$2n$ is detected in $\HAQ_{2n}(MU,S;\F_2)$ if and only if
$n$ is a power of~$2$. Therefore $MU_{(2)}$ is not minimal atomic.

For an odd prime $p$, there is splitting of infinite loop spaces
\[
BU_{(p)} \simeq W_1 \times\cdots\times W_{p-1},
\]
where $W_r$ is the zeroth space of $\Sigma^{2r}\ell$, so that the
natural map
\[
\Sigma^\infty BU_{(p)} \lra \Sigma^2ku_{(p)}
              \simeq \bigvee_{1\leq r\leq p-1} \Sigma^{2r}\ell
\]
factors through a wedge of maps $\Sigma^\infty W_r\lra\Sigma^{2r}\ell$.
An exercise in using the definition of the Adams splitting together with
the fact that the lowest degree cohomology class for $W_r$ corresponds to
the Newton polynomial in the Chern classes in $H^{2r}(BU;\F_p)$, shows
that the induced composition
\[
\Sigma^\infty \CP^\infty \xrightarrow{\;x^r\;} \Sigma^\infty W_r
                                           \lra \Sigma^{2r}H\F_p
\]
has the effect on homology of sending $\beta_n$ to the coefficient
of $t^n$ in the series
\[
\xi(t)^r = \bigl(\sum_{s\geq0}\xi_st^{p^s}\bigr)^r.
\]
In particular, this shows that when $1\leq r \leq p-1$, each generator
$b_{p^s-1+r}\in H_{2(p^s-1+r)}(MU;\F_p)$ survives to give a non-zero
element in $H_{2(p^s-1+r)}(\Sigma^2ku;\F_p)$, therefore the corresponding
generator $x_{p^s-1+r}\in\pi_{2(p^s-1+r)}MU$ also survives to a non-zero
element of $H_{2(p^s-1+r)}(\Sigma^2 ku;\F_p)$.

To summarise, the Hurewicz homomorphism $\theta\:\pi_*MU\lra\HAQ_*(MU,S;\F_p)$
detects homotopy, so $MU_{(p)}$ is not minimal atomic for any prime~$p$.

\bigskip
\noindent
\textbf{A core for $MU_{(2)}$:}
The problem of identifying a core for $MU_{(p)}$ was first studied
in~\cite{Hu-Kriz-May} where the following example for the case $p=2$
appeared. We give a new verification that the domain is indeed minimal
atomic. Of course, it would be very interesting to identify cores of
$MU$ for the odd primes.

Consider the infinite loop map $BU\lra BSp$ classifying quaternionification
of complex bundles. The fibre is also an infinite loop map $j\:Sp/U\lra BU$,
where $Sp/U$ is the zeroth space of $\Sigma^2 ko$. Then $j$ induces a Thom
spectrum $Mj$ which is a commutative $S$-algebra and the natural map
$Mj\lra MU$ is a morphism of commutative $S$-algebras. Furthermore, the
associated map in $2$-local homology gives an isomorphism onto half of
$H_*(MU;\Z_{(2)})$:
\[
H_*(Mj;\Z_{(2)}) = \Z_{(2)}[y_{2r-1}:r\geq1] \xrightarrow{\iso}
\Z_{(2)}[y'_{2r-1}:r\geq1] \subseteq H_*(MU;\Z_{(2)}),
\]
where $y_{2r-1},y'_{2r-1}$ have degree $2r-2$ and in $H_*(MU;\Z_{(2)})$,
\[
y'_{2r-1} \equiv b_{2r-1} \pmod{\mathrm{decomposables}}.
\]
There is a similar result for mod~$2$ homology.

{}From now on in this example, all spectra are assumed to be
localised at~$2$ and we drop this from the notation.

The argument of~\cite{Hu-Kriz-May} shows that $Mj$ is a wedge
of suspensions of $BP$ and as the induced map in homology is
a monomorphism so is that in homotopy. So to show that
$Mj\lra MU$ is a core we only need to show that $Mj$ is minimal
atomic.

This time we need to examine the mod~$2$ $\HAQ$ Hurewicz
homomorphism which amounts to a homomorphism
\[
\pi_*Mj \lra H_*(\Sigma^2 ko;\F_2).
\]
By work of Stong~\cite{Stong:Book}, the map induced by the bottom
cohomology class gives a homology monomorphism
$H_*(\Sigma^2 ko;\F_2)\lra \mathcal{A}(2)_*$. There is a commutative
diagram
\[
\begin{CD}
\Sigma^\infty Sp/U @>>> \Sigma^\infty BU \\
@VVV  @VVV \\
\Sigma^2 ko @>>> \Sigma^2 ku
\end{CD}
\]
and on applying $H_*(-;\F_2)$ and using the calculations described
for~$MU$, we see that
\[
\sigma_*(y'_{2r-1}) = 0\quad(r\geq1).
\]
Hence the mod~$2$ $\HAQ$ Hurewicz homomorphism is also trivial in
positive degrees.

We remark that the morphisms of $S$-algebras $S\lra Mj\lra MU$ give
rise to a cofibre sequence of $MU$-modules
\[
MU\wedge\Sigma^2 ko \lra MU\wedge\Sigma^2 ku \lra \Omega_{MU/Mj},
\]
and this is equivalent to the cofibre sequence
\[
MU\wedge\Sigma^2 ko \lra MU\wedge\Sigma^2 ku \lra MU\wedge\Sigma^4 ko,
\]
where we view $ku$ as $ko\wedge\cone_\eta$ for $\eta\in\pi_1S$ the
generator. See~\cite{JR:Galois}*{proposition~5.3.1} for more on this
sequence in a Galois theoretic context. Smashing with $H\F_2$ over
$MU$ now gives the usual short exact sequence
\[
0\ra H_*(ko;\F_2)\lra H_*(ku;\F_2)\lra H_{*-2}(ko;\F_2) \ra0.
\]

\bigskip
\noindent
\textbf{$MSp$ $2$-locally:}
Here only the prime $p=2$ presents an interesting question.
\begin{proposition}\label{prop:MSp}
$MSp_{(2)}$ is minimal atomic as a commutative $S$-algebra
but not as an $S$-module.
\end{proposition}
\begin{proof}
We will show that the mod~$2$ $\HAQ$ Hurewicz homomorphism
is trivial in positive degrees. It is known that the mod~$2$
Hurewicz homomorphism $\pi_*MSp \lra H_*(MSp;\F_2)$ is not
trivial in positive degrees, so $MSp$ is not minimal atomic
as a $2$-local spectrum.

For this example, $BSp$ is the zeroth space of the $3$-connected
spectrum $\Sigma^4 ko$. From~\cite{Stong:Book}, the bottom integral
cohomology class induces a monomorphism
\[
H_*(ksp;\F_2)\lra H_*(\Sigma^4 H\F_2;\F_2) = \mathcal{A}(2)_{*-4}.
\]
However, the crux of our argument involves a beautiful result
of Floyd~\cite{Floyd} and we describe this in some detail.

Recall that the natural map $MSp \lra MO$ induces a mod~$2$
homology isomorphism onto the $4$-th powers:
\begin{equation}\label{eqn:H*MSp->H*MO}
H_*(MSp;\F_2) \xrightarrow{\iso} H_*(MO;\F_2)^{(4)}
                                     \subseteq H_*(MO;\F_2).
\end{equation}
In~\cite{Floyd}*{theorems~5.3,5.5}, it is shown that $\pi_*MO$
has a family of polynomial generators $z_r\in\pi_r MO$ ($r+1$
not a power of~$2$) for which the polynomial subring
\[
P_* = \F_2[z_r^{\kappa(r)}:\text{$r+1$ not a power of~$2$}]
                                           \subseteq \pi_*MO,
\]
satisfies
\[
\im (\pi_*MSp \lra \pi_*MO) \subseteq P_*^{(8)},
\]
where
\[
\kappa(r) =
\begin{cases}
2 & \text{if $r$ is odd or a power of $2$}, \\
1 & \text{if $r$ is even and not a power of $2$}.
\end{cases}
\]
Kochman~\cite{Kochman:I} shows that these two rings are actually
equal, but for our purposes it suffices to have the inclusion.

For our purposes we need only remark that~\eqref{eqn:H*MSp->H*MO}
together with the fact that the Hurewicz homomorphism
$\pi_*MO\lra H_*(MO;\F_2)$ is monic, combine to show that the
Hurewicz homomorphism for $MSp$ has image contained in
the squares:
\[
\im (\pi_*MSp \lra H_*(MSp;\F_2)) \subseteq H_*(MSp;\F_2)^{(2)},
\]
hence in positive degree this image is contained in the
decomposables. Therefore the composition
\[
H_*(MSp;\F_2) \xrightarrow{\iso} H_*(BSp;\F_2)
\xrightarrow{\sigma_*} H_*(\Sigma^4 ksp;\F_2)
\]
is trivial in positive degrees since the second map annihilates
decomposables.
\end{proof}

\section*{Concluding remarks}

An important goal for this paper was to show how $\TAQ$ could
be used to introduce cellular arguments into the study of
commutative $S$-algebras, and in particular to aid in identifying
the minimal atomic $p$-local $S$-algebras introduced
in~\cite{Hu-Kriz-May} using ideas suggested by~\cite{AJB-JPM}.
We have simplified things by confining attention to simply
connected examples, however the theory can be adapted to work
as well when $\pi_0R$ is finite cyclic of prime power order,
and we have already investigated some examples. It would also
be interesting to have an analogous theory for cellular
$S$-algebras in Bousfield localised situations, for example
for $K(n)$-local algebras. We expect to address these issues
in future work.

\end{document}